\documentclass[11pt]{amsart}
\usepackage{comment}
\usepackage{amssymb}
\usepackage{graphicx}
\usepackage{mathtools}
\usepackage{subcaption}
\usepackage{tikz-cd}
\usepackage{algorithm2e}

\let\visiblecomments y

\def\articletheorems{
\newtheorem{thm}{Theorem}[section]
\newtheorem{lem}[thm]{Lemma}

\newtheorem{defn}[thm]{Definition}
\newtheorem{cor}[thm]{Corollary}
\newtheorem{prop}[thm]{Proposition}
\newtheorem{ex}[thm]{Example}

}
\def\setof#1{\mbox{$\{\,#1\,\}$}}

\def\mathobj#1{\mbox{$#1$}}
\def\NN{\mathobj{\mathbb{N}}}

\def\ZZ{\mathobj{\mathbb{Z}}}
\def\PP{\mathobj{\mathbb{P}}}

\def\cP{\text{$\mathcal P$}}

\def\cM{\text{$\mathcal M$}}

\newcommand{\cl}{\operatorname{cl}}
\newcommand{\bd}{\operatorname{bd}}
\newcommand{\Inv}{\operatorname{Inv}}
\newcommand{\Int}{\operatorname{int}}

\newcommand{\im}{\operatorname{im}}
\newcommand{\dom}{\operatorname{dom}}

\newcommand{\mto}{\multimap}

\articletheorems

\usepackage[cp1250]{inputenc}

\numberwithin{table}{section}

\begin{document}
\title[Morse equation for discrete multivalued dynamical systems]{The Morse equation in the Conley index theory for discrete multivalued dynamical systems}
\author[B. Batko]{Bogdan Batko}
\address{Bogdan Batko,
 Department of Mathematics,
Rutgers University,
Piscataway, NJ 08854, USA}
\email{bogdan.batko@rutgers.edu}
\address{{\it On leave from}}
\address{Division of Computational Mathematics,
 Faculty of Mathematics and Computer Science,
         Jagiellonian University,
         ul.~St. \L{}ojasiewicza 6, 30-348~Krak\'ow, Poland}
\email{bogdan.batko@uj.edu.pl}

\thanks{
   The author was partially supported by the Polish National Science Center under Ma\-estro Grant No. 2014/14/A/ST1/00453, Opus Grant No. 2019/35/B/ST1/00874, and the DARPA contract HR0011-16-2-0033.
}

\begin{abstract} 
A recent generalization of the Conley index to discrete multivalued dynamical systems without a continuous selector is motivated by applications to data--driven dynamics. In the present paper we continue the program by studying attractor--repeller pairs and Morse decompositions in this setting. In particular, we prove Morse equation and Morse inequalities.  
\end{abstract}
\maketitle
\section{Introduction}
The present data--driven world sets new challenges to contemporary science. There is a growing interest in coarse theories capable to extract robust information hidden in noisy experimental data. A good example is the rapid development in persistent homology (cf. \cite{EH2010}) invented to investigate homological features of topological spaces known only from a cloud of points or samples. Conley theory (cf. \cite{C78}) provides a robust topological invariant for studying dynamical systems. It has been used to prove the existence of stationary and periodic solutions, heteroclinic connections and chaotic invariant sets. The potential of Conley theory in the context of data became apparent via a generalization to multivalued dynamical systems \cite{KM95} proposed as a tool in computer assisted proofs in dynamics. Multivalued dynamics,  important on its own right (cf. e.g. \cite{AC84}), becomes in such proofs a tool in the study of single valued dynamics (cf. e.g. \cite{MiMr95,MR95,TW2002,MZ2000,KMM}). However, the practical use of the generalization \cite{KM95} has been severely restricted by some strong assumptions difficult to be fulfilled in practice. Fortunately, a recent revision of the theory (cf. \cite{BM2016,B2017}) removes these limitations opening the way to applications in data--driven dynamics.    

Assume time series data have been collected as a result of a measurement, an experiment, or an observation of an unknown dynamical system. Using similar techniques as in \cite{MiMrReSz99,Szymczak-1997,BMMP2020,BGHKMV} one can construct a multivalued map that represents an unknown generator of the underlying system. Both, toy examples (cf. \cite{BM2016,B2017}) and constructions for real applications (cf. \cite{BMMP2020,BMMP2020app}) show that usually such a map does not admit a continuous selector. Nevertheless, we can identify isolated invariant sets, and hire the Conley index theory for multivalued maps, because the construction of the index works under minimal assumptions on a multivalued map; in particular, it does not require single valued selectors. Eventually, we extend the obtained results to the underlying unknown single valued dynamical system using continuation properties of the index. 

In the present paper we continue the program that we started in \cite{BM2016,B2017}. Heading towards the comprehensive description of the dynamics we need to gain more insight into the internal structure of isolated invariant sets. One of the relevant descriptors is a Morse decomposition of an isolated invariant set and the associated Morse equation. Leaving mathematical sophistication aside, one can think that the Morse decomposition is a decomposition of a given isolated invariant set into a finite union of pairwise disjoint isolated invariant sets, called Morse sets, and connecting orbits between them, such that outside the Morse sets the dynamics is gradient--like. 
There is a rich variety of approaches to global dynamics that involve
Morse decompositions; see for instance \cite{RFr86,RFrKMi88,RFr89,BD99} or \cite{BK88,Li2007} for continuous-time dynamical
systems without uniqueness. For recent results that combine classical and combinatorial
dynamics we refer to \cite{MM2017,DJKKLM,BKMW2020}. An approach to the poset structure of Morse decomposition
with the emphasis on lattice structures of attractors, for both an underlying dynamical system
and its combinatorialization - a combinatorial multivalued map, is demonstrated in the series
of papers \cite{KMV05, KMV1, KMV2, KMV3}. The presented list of papers is not complete and does not pretend to be
complete in any sense.

The Morse equation describes the relationship between the Conley index of an isolated invariant set and the Conley indices of its Morse sets. In particular, information on (co)homologically nontrivial connections between Morse sets can be derived.
The classical Morse inequalities concern nondegenerate critical points of a gradient flow on a compact manifold, and show the correspondence between the $k$-th Betti number of the manifold and the number of critical points with Morse index $k$, that is critical points with the $k$-dimensional unstable invariant manifold (cf. e.g. \cite{RB80,RB82}).  One of the possible generalizations of the classical inequalities due to Morse and to Smale (cf. \cite{Sm60}) in the Conley index theory was developed by Conley and Zehnder to flows (cf. \cite{CZ84}). Afterwords, the Morse equation in the Conley index theory was proved by Rybakowski in \cite{RZ85} for semiflows, and  in the discrete time case by Franks \cite{JFr82} and Mrozek \cite{MMr91}. 	  

Assume that we are given a (continuous-- or discrete--time) dynamical system on a locally compact metric space. A set $S$ is said to be isolated invariant if it is the maximal invariant set contained in some compact neighborhood of itself. In the cohomological Conley index theory with such a set one associates a special pair of sets, called an index pair. Then, the cohomological Conley index of $S$ is defined to be the Alexander--Spanier cohomology of the index pair. One proves that, in line with the need, this is an invariant of $S$. Thus, we can associate with $S$ the Poincar\'e series $p(t,S)$, the power series in $t$ with the ranks of the cohomology modules as coefficients. Now, assume that $\cM:=\{M_i\ |\ i\in\{1,2,\dots,n\}\}$ is a Morse decomposition of $S$. Hence, in particular, the sets $M_i$ are pairwise disjoint isolated and invariant subsets of $S$. Then, the Morse equation takes the form
$$
\sum _{i=1} ^n p(t,M_i)=p(t,S)+(1+t)Q(t),
$$ 
where $Q$ is a formal power series with nonnegative integer coefficients (cf. \cite{CZ84}). The terms in $Q$ provide information about nontrivial connections between pairs of Morse sets. One can observe that the above equation generalizes the classical Morse inequalities
(cf. e.g. \cite{RZ85,RS99}).

The aim of this paper is to prove the Morse equation in the Conley index theory for discrete multivalued dynamical systems.

The organization of the paper is as follows. In Section \ref{sec:preliminary} we provide basic definitions related to the Conley index theory for multivalued maps. Section \ref{sec:A-R} is devoted to Morse decompositions. We define repeller--attractor pairs and provide a characterization of Morse sets via the associated sequence of attractors. The key step for the prove of the Morse equation is the construction of the, so called, index triple for a given repeller--attractor pair. This is presented in Section \ref{sec:triples}. Finally, we prove the Morse equation in Section \ref{sec:Morse_eq}.   

\section{Preliminaries}
\label{sec:preliminary}
Throughout the paper $\ZZ$ and $\NN$ will stand for the sets of all integers and all positive integers, respectively. By an interval in $\ZZ$ we mean the trace of a real interval in $\ZZ$.

Given a topological space $X$ and a subset $A\subset X$, by $\Int_X A$, $\cl_X A$ we denote the {\em interior} of $A$ in $X$ and the {\em closure} of $A$ in $X$ respectively. We omit the subscript $X$ if the space is clear from the context. 

Let $X$, $Y$ be topological spaces, and let $F\colon X\mto Y$ denote a multivalued map, that is a map $F\colon X\ni x\mapsto F(x)\in \cP(Y)$, where $\cP(Y)$ is the power set of $Y$. We define an {\em effective domain} of $F$ by 
$
\dom (F):=\{x\in X\ |\ F(x)\neq\emptyset\}.
$
A multivalued map $F$ is {\em upper semicontinuous} if for any closed $B\subset Y$ its large counter image under $F$, that is the set $F^{-1}(B):=\{x\in X\ |\ F(x)\cap B\neq\emptyset\}$, is closed. 

Assume that a multivalued self-map $F\colon X\mto X$ is given.

Let $I$ be an interval in $\ZZ$ with $0\in I$. We say that a single valued mapping $\sigma\colon I\to X$ is a \emph{solution for $F$ through $x\in X$} if $\sigma(n+1)\in F(\sigma(n))$ for all $n, n+1\in I$ and $\sigma(0)=x$.
Given a subset $N\subset X$, the set 
$$
\Inv(N,F):=\{x\in N\ |\ \exists\sigma\colon\ZZ\to N \text{ a solution for } F \text{ through } x\}
$$ 
is called the {\em invariant part} of $N$. A compact subset $N\subset X$ is an {\em isolating neighborhood} for $F$ if $\Inv(N,F)\subset\Int N$. A compact set $S\subset X$ is said to be {\em invariant} with respect to $F$ if $S = \Inv(S,F)$. It is called an {\em isolated invariant set} if it admits an isolating neighborhood $N$ for $F$ such that $S = \Inv(N,F)$ (cf. \cite[Definition 4.1, Definition 4.3]{BM2016}). 
\begin{defn}[cf. {\cite[Definition 4.7]{BM2016}}]
\label{defn:wip}
Let $N\subset X$ be an isolating neighborhood for $F$.  A pair $P=(P_1,P_2)$ of compact sets $P_2\subset P_1\subset N$ is called a {\em weak index pair} in $N$ if
\begin{itemize}
\item[(a)] $F(P_i)\cap N\subset P_i$ for $i\in\{1,2\}$,
\item[(b)] $\bd_F P_1:= \cl P_1 \cap \cl(F(P_1)\setminus P_1)\subset P_2$,
\item[(c)] $\Inv(N,F)\subset \Int(P_1\setminus P_2)$,
\item[(d)] $P_1\setminus P_2\subset \Int N$.
\end{itemize}
\end{defn}

Given a weak index pair $P$ in an isolating neighborhood $N\subset X$ for $F$ we set
$$
T_N(P):=(T_{N,1}(P),T_{N,2}(P)):=(P_1\cup (X\setminus\Int N),P_2\cup (X\setminus\Int N)).
$$
Assume $F$ is determined by a morphism, which in particular holds true if $F$ has acyclic values (cf. \cite{G76}). Following \cite{BM2016} recall that $F_P$, the restriction of $F$ to the domain $P$, is a multivalued map of pairs, $F_P\colon P\mto T_N(P)$; the inclusion $i_P:P\to T_N(P)$ induces an isomorphism in the Alexander--Spanier cohomology; and the {\em index map} $I_{F_P}$ is defined as an endomorphism of $H^*(P)$ given by
$$
I_{F_P}=F_P ^*\circ (i_P ^*)^{-1}.
$$
Applying the Leray functor $L$ (cf. \cite{M90}) to the pair $(H^*(P),I_{F_P})$ we obtain a graded module $L(H^*(P),I_{F_P})$ with its endomorphism, called the {\em Leray reduction of the Alexander--Spanier cohomology of a weak index pair} $P$. By \cite[Definition 6.3]{BM2016}, this is the {\em cohomological Conley index of} $\Inv(N,F)$, which we denote by $C(\Inv(N,F),F)$.

We introduce the following notation for future use. Given $N\subset X$, $x\in N$ and $n\in\ZZ^+$ we put:
\begin{eqnarray*}
F_{N,n}(x)&:=& \{y\in N\mid \exists\,\sigma: [0,n]\to N\mbox{ a solution for }F\\
&&\mbox{\hspace{4.25cm}with }\sigma(0)=x,\,\sigma(n)=y\},\\
F_{N,-n}(x)&:=& \{y\in N\mid \exists\,\sigma: [-n,0]\to N\mbox{ a solution for }F\\&&\mbox{\hspace{4.55cm}with }\sigma(-n)=y,\,\sigma(0)=x\},\\
F_N^+(x)&:=&\bigcup\setof{F_{N,n}(x)\,|\,n\in\ZZ^+},\\
F_N ^-(x)&:=&\bigcup\setof{F_{N,-n}(x)\,|\, n\in\ZZ^+}.
\end{eqnarray*}

\section{Morse decompositions and repeller--attractor pairs}\label{sec:A-R}
Given a solution $\sigma:\ZZ\to X$ of a multivalued upper semicontinuous map $F:X \mto X$,
we define its {\em $\alpha$-- and $\omega$--limit sets} respectively by
$$
\alpha(\sigma):=\bigcap_{k\in\ZZ}\cl\sigma((-\infty,k]),  \quad\quad  \omega(\sigma):=\bigcap_{k\in\ZZ}\cl\sigma([k,+\infty)).
$$
Let us point out that unlike in the single valued case, we define $\alpha$-- and $\omega$--limit sets for a given solution $\sigma$ through an $x\in X$, not for an $x$ itself, because backward nor  forward solutions through $x$ need not be unique.
\begin{defn}[cf. {\cite[Definition 3.9]{BKMW2020}}]\label{def:Morse_dec}
{\rm
Let $S$ be an isolated invariant set of $F:X\mto X$.
We say that the family $\cM:=\{M_r\ |\ r\in\PP\}$ indexed by a poset $(\PP,\leq)$
is a {\em Morse decomposition of $S$} if  the following conditions are satisfied:
\begin{enumerate}
\item the elements of $\cM$ are mutually disjoint isolated invariant subsets of $S$,
\item for every full solution $\sigma$ in $X$ there exist $r,r'\in\PP$, $r\leq r'$, such that $\alpha(\sigma)\subset M_{r'}$ and $\omega(\sigma)\subset M_r$,
\item if for a full solution $\sigma$ in $X$ and $r\in\PP$ we have $\alpha(\sigma)\cup\omega(\sigma)\subset M_r$, then $\im\sigma\subset M_r$.
\end{enumerate}
}
\end{defn}
\begin{figure}[]
\begin{minipage}[t]{0.47\linewidth}
\centering
\includegraphics[width=\textwidth]{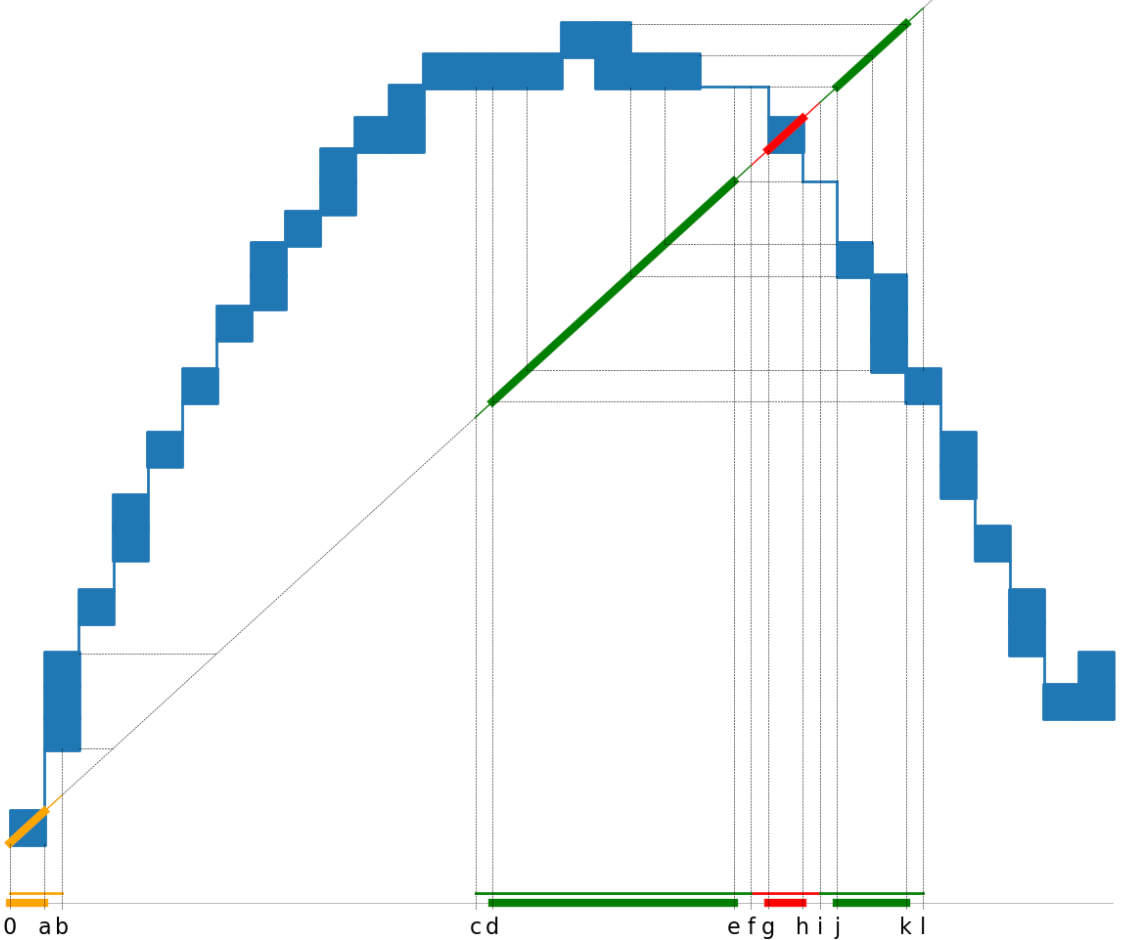}
\subcaption{Morse decomposition $\cM=\{M_1,M_2,M_3\}$, with $M_3=[0,a]$, $M_2=[g,h]$, and $M_1=[d,e]\cup[j,k]$, of an isolated invariant set $S=[0,k]$ of $F$. Morse sets depicted as thick lines in orange, red and green, respectively, and the associated isolating neighborhoods as thin lines in the same colors.} 
\end{minipage}
\hspace{4mm}
\begin{minipage}[t]{0.47\linewidth}
\centering
\includegraphics[width=\textwidth]{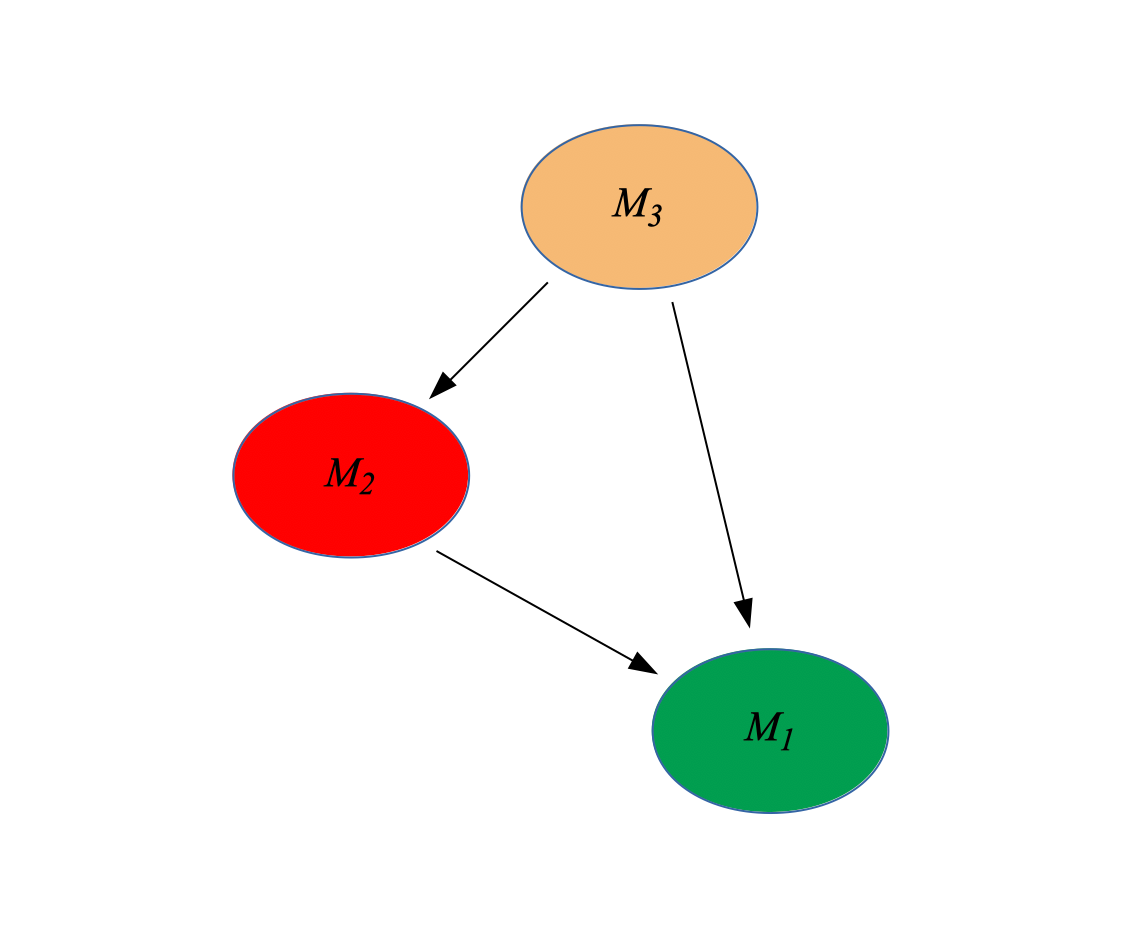}
\subcaption{Morse graph of $F$ - a directed graph with the set of vertices $\cM$ and arrows presenting the existence of connecting orbits.}
\end{minipage}\\
\begin{minipage}[t]{0.47\linewidth}
\centering
\includegraphics[width=\textwidth]{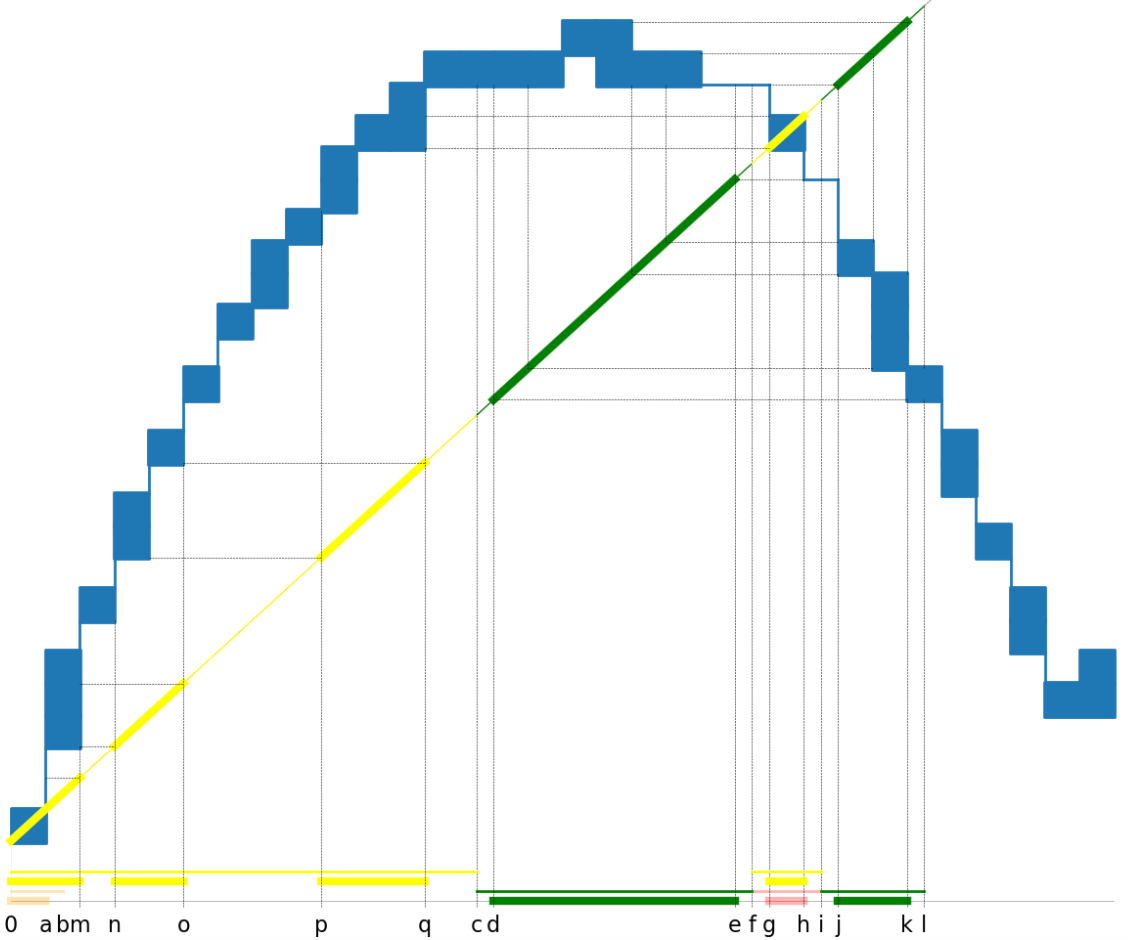}
\subcaption{Repeller--attractor pair $(A_1 ^*,A_1)$ in $S$. Attractor $A_1=M_1=[d,e]\cup [j,k]$ and its trapping region $[c,f]\cup [i,l]$ depicted as green lines. The dual repeller $A_1 ^*=[0,m]\cup[n,o]\cup [p,q]\cup[g,h]$ and its isolating neighborhood in yellow.}  
\end{minipage}
\hspace{5mm}
\begin{minipage}[t]{0.47\linewidth}
\centering
\includegraphics[width=\textwidth]{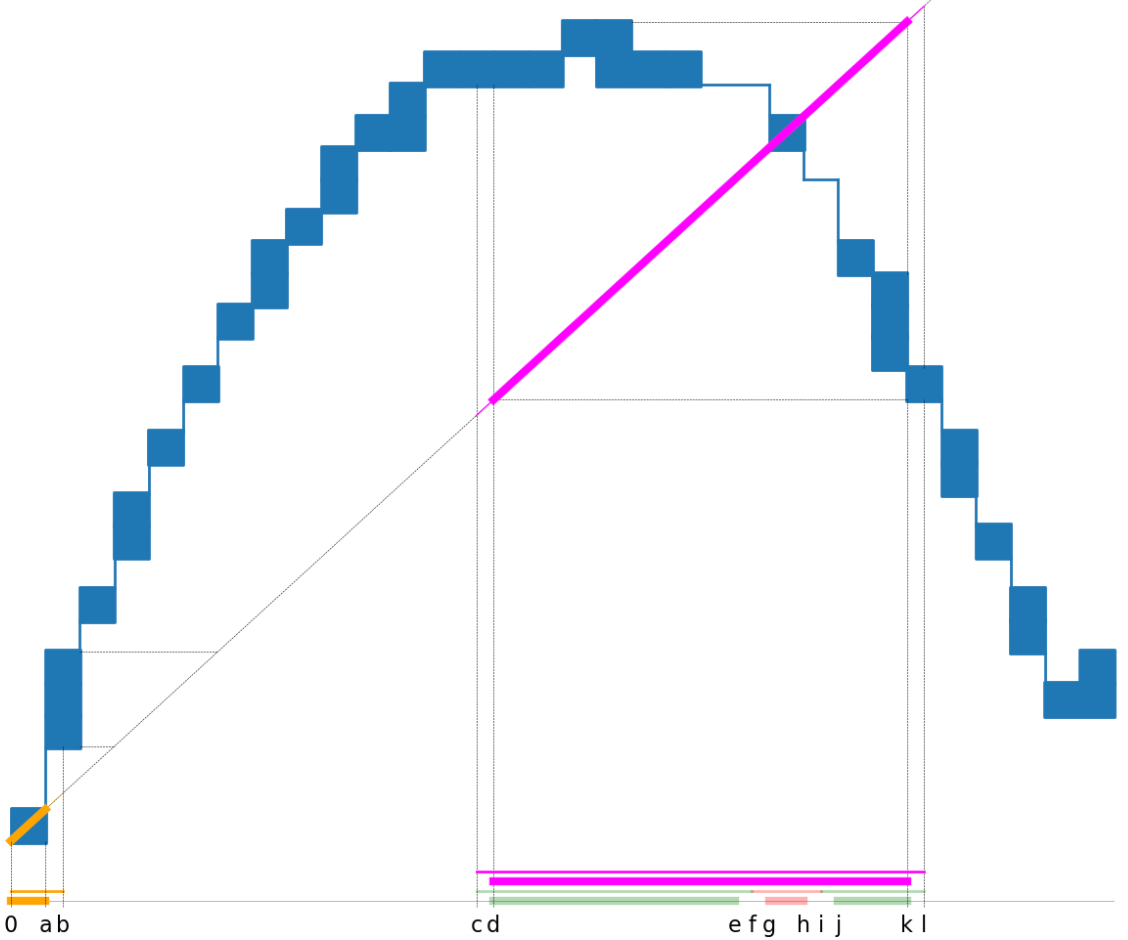}
\subcaption{Repeller--attractor pair $(A_2 ^*,A_2)$ in $S$. Attractor $A_2=[d,k]$ and its trapping region $[c,l]$ in magenta. The dual repeller $A_2 ^*=[0,a]=M_3$ and its isolating neighborhood depicted as orange lines. Note that $[0,c]$ is an another isolating neighborhood of $A_2 ^*$ (cf. Lemma \ref{lem:rep-attr}).}
\end{minipage}
\caption{Graph of a multivalued map $F:[0,1]\mto [0,1]$ in blue. Panel (A) - Morse decomposition of an isolated invariant set $S=[0,k]$ for $F$, panel (B) - Morse graph, panel (C) and panel (D) -  repeller--attractor pairs in $S$.}
\label{fig:Morse_dec}
\end{figure}
The partial order $\leq$ on $\PP$ will be called an {\em admissible} ordering of the Morse decomposition $\cM$. Note that it is not uniquely determined. Moreover, there is an "extremal" admissible ordering $\leq_F$, given by $p\leq_F q$ if and only if there exists a sequence of distinct elements $p=r_0,r_1,\dots,r_k=q$ of $\PP$ such that, for any $j\in\{1,2,\dots,k\}$, there exists a solution $\sigma$ for $F$ with $\alpha(\sigma)\subset M_{r_j}$ and $\omega(\sigma)\subset M_{r_{j-1}}$. The ordering $\leq_F$ is extremal in the sense that any admissible ordering of $\cM$ is an extension of $\leq_F$. One can observe that for any admissible ordering $\leq$ there exists its linear extension which is also admissible. In the latter case, that is whenever the admissible ordering is linear, for simplicity we write $\cM=\{M_1,M_2,\dots, M_n\}$.    
\begin{ex}\label{ex:F}
{\rm Consider multivalued map $F:[0,1]\mto[0,1]$, graph of which is presented in Figure \ref{fig:Morse_dec}. Note that $S=[0,k]$ is an isolated invariant set for $F$ and $\cM=\{M_1,M_2,M_3\}$, with $M_1=[0,a]$, $M_2=[g,h]$, and $M_3=[d,e]\cup[j,k]$, is a Morse decomposition of $S$ (see Figure \ref{fig:Morse_dec}(A) and (B)).}
\end{ex}
\begin{defn}\label{def:attractor} 
{\rm We say that an isolating neighborhood $T$ for $F$ is a {\em trapping region} if $F(T)\subset T$. A subset $A$ of an isolated invariant set $S$ is called an {\em attractor} (in $S$) if it admits a trapping region $T$ which isolates $A$ (relative to $S$). Given an attractor $A$ in $S$, the set $A^*:=\{x\in S\ |\mbox{ there exists a solution }\sigma:\ZZ\to S\mbox{ through }x\mbox{ with }\omega(\sigma)\cap A=\emptyset\}$ will be called the {\em repeller} dual to $A$
, and the pair $(A^*,A)$ will be called the {\em repeller--attractor pair} in $S$.} 
\end{defn}
\begin{ex}\label{ex:seq_attr}
{\rm Consider Morse decomposition of an isolated invariant set $S=[0,k]$ with respect to a multivalued map presented in Example \ref{ex:F} (see Figure \ref{fig:Morse_dec}). Note that the pairs of sets $(A_0 ^*,A_0)=(S,\emptyset)$, $(A_1 ^*, A_1)=([0,m]\cup[n,o]\cup [p,q]\cup[g,h],[d,e]\cup[j,k])$ (Figure \ref{fig:Morse_dec}(C)), $(A_2 ^*,A_2)=([0,a],[d,k])$ (Figure \ref{fig:Morse_dec}(D)) and $(A_3 ^*,A_3)=(\emptyset, S)$ are repeller-attractor pairs in $S$. 

Observe that $\cM$, the sequence of attractors $\emptyset=A_1\subset A_2\subset A_3=S$,  and the associated sequence of dual repellers $S=A_0 ^*\supset A_1 ^*\supset A_2 ^*\supset A_3 ^*=\emptyset$, satisfy all the assertions of upcoming Theorem \ref{thm:Morse_attr-rep}, Proposition \ref{prop:repattr-in-a} and Theorem \ref{thm:Attr-rep-Morse}.
}
\end{ex}

Assume that $A$ is an attractor and $T$ is a given trapping region of $A$. It is easily seen that if $\sigma:\ZZ\to S$ is a  solution  for $F$ through an $x\in T$, say $\sigma(0)=x$, then $\sigma(\ZZ^+)\subset T$. Consequently, $\omega(\sigma)\subset T$ according to the compactness of $T$. It follows that  $\omega(\sigma) \subset A$, because $\omega(\sigma)$ is invariant and $A=\Inv(T,F)$. We state this simple observation as a lemma for further references.

\begin{lem}\label{lem:omega_in_A}
Let $A$ be an attractor in $S$ with a trapping region $T$. If $\sigma:\ZZ\to S$ is a solution for $F$ with $\im\sigma\cap T\neq\emptyset$ then $\omega(\sigma)\subset A$.
\end{lem}

\begin{lem}\label{lem:rep-attr}
Let $(A^*,A)$ be a repeller--attractor pair in $S$. Then $A^*$ and $A$ are disjoint isolated invariant sets. Moreover, if $T$ is a trapping region for $A$ then $N:=S\setminus\Int T$ is an isolating neighborhood for $A^*$.
\end{lem}
\proof
Directly by the definition, $A$ is an isolated invariant set, and $A^*\cap A=\emptyset$. It is easy to see that $A^*\cap T=\emptyset$. Indeed, otherwise we would have an $x\in A^*\cap T$, and a solution $\sigma:\ZZ\to S$ with $\sigma(0)=x$ with $\omega(\sigma)\cap A=\emptyset$, as $x\in A^*$. However, $x\in\im\sigma\cap T$, which according to Lemma \ref{lem:omega_in_A} means that $\omega(\sigma) \subset A$, a contradiction. 

We shall prove that $N$ is an isolating neighborhood for $A^*$. Clearly, $N$ is compact. Let $x\in A^*$. There exists a solution $\sigma:\ZZ\to S$ through $x$ with $\omega(\sigma)\cap A=\emptyset$. Note that $\im\sigma\subset\Int  N$, because otherwise we would have $\im\sigma\cap T\neq\emptyset$ and in turn, by Lemma \ref{lem:omega_in_A}, $\omega(\sigma)\subset A$, a contradiction. This shows the inclusions $A^*\subset \Inv N$ and $A^*\subset \Int N$. Now consider an $x\in \Inv N$ and a solution $\sigma:\ZZ\to N$ through $x$. By compactness of $N$ we get $\omega(\sigma)\subset N$, showing that $\omega(\sigma)\cap A=\emptyset$. Thus, $x\in A^*$, showing the inclusion $\Inv N\subset A^*$. This completes the proof.
\qed
\begin{lem}\label{lem:rep-attr-sol}
Let $(A^*,A)$ be a repeller--attractor pair in $S$, let $\sigma:\ZZ\to S$ be a solution for $F$ through an $x\in S$, and  let $T$ be a trapping region for $A$. The following holds true:
\begin{itemize}
\item[(i)] if $x\notin T$ then $\alpha(\sigma)\subset A^*$,
\item[(ii)] if $x\notin A^*\cup A$ then $\alpha(\sigma)\subset A^*$ and $\omega(\sigma)\subset A$,
\item[(iii)] if $\omega(\sigma)\cap A^*\neq\emptyset$ then $\im\sigma\subset A^*$,
\item[(iv)] if $\alpha(\sigma)\cap A\neq\emptyset$ then $\im\sigma\subset A$.
\end{itemize}
\end{lem}
\proof
Set $N:=S\setminus\Int T$ and, without loss of generality, assume that $\sigma(0)=x$.

For the proof of (i), observe that $\sigma(\ZZ^-)\subset N$. By the compactness of $N$ we infer that $\alpha(\sigma)\subset N$. Consequently, $\alpha(\sigma) \subset A^*$, because $\alpha(\sigma)$ is invariant and, according to Lemma \ref{lem:rep-attr}, $\Inv(N,F)=A^*$. 

We shall verify (ii). Since $x\notin A^*$, we get $\omega (\sigma)\cap A\neq\emptyset$. Then $\im\sigma\cap T\neq\emptyset$, and by Lemma \ref{lem:omega_in_A}, $\omega(\sigma)\subset A$. If $\im\sigma\not\subset T$ then the inclusion $\alpha(\sigma)\subset A^*$ follows from (i). We have left  the case $\im\sigma\subset T$. However, in that case we would have $\im\sigma\subset A$, and in turn $x\in A$, leading to a contradiction.

We verify (iii). For contradiction suppose that $\im\sigma\not\subset A^*$. By Lemma \ref{lem:rep-attr}, $N:=S\setminus\Int T$ is an isolating neighborhood for $A^*$, therefore we have $\im\sigma\not\subset N$. This means that $\im\sigma\cap T\neq\emptyset$. However, by Lemma \ref{lem:omega_in_A} we get $\omega(\sigma)\subset A$, a contradiction.

For the proof of (iv) observe that the assertion holds whenever $\im\sigma\subset T$. Thus, assume the contrary. But then, by (i), we get $\alpha(\sigma)\subset A^*$, a contradiction.
\qed
\begin{lem}\label{lem:rep_neighborhood}
Let $K$ be an isolating neighborhood of $S$, let $(A^*,A)$ be a repeller--attractor pair in $S$, and let $N\subset K$ be a compact neighborhood of $A^*$ disjoint from $A$. Then $N$ is an isolating neighborhood for $A^*$.
\end{lem}
\proof 
Since $A^*\subset\Int N$ and $N$ is compact, it suffices to verify that $A^*=\Inv (N,F)$.

Inclusion $A^*\subset\Inv (N,F)$ is straightforward, as $A^*$ is invariant and $A^*\subset N$. To see the opposite inclusion consider an $x\in \Inv (N,F)$ and a solution $\sigma$ through $x$ in $N\subset K$. Then, by compactness of $N$ we have $\omega(\sigma)\subset N$. In particular, $\omega(\sigma)\cap A=\emptyset$, as $N$ and $A$ are disjoint. This along with Lemma \ref{lem:rep-attr-sol} yields $x\in A^*$, and completes the proof.
\qed
\begin{thm}\label{thm:Morse_attr-rep}
Let $\cM:=\{M_1,M_2,\dots,M_n\}$ 
be a Morse decomposition of an isolated invariant set $S$ with respect to an upper semicontinuous $F:X\mto X$. Then there exists a family 
$\emptyset =A_0\subset A_1\subset A_2\subset\dots\subset A_n=S$ of attractors in $S$ such that $M_j=A_j\cap A^* _{j-1}$ for $j\in\{1,2,\dots,n\}$, where $A^* _j$ is a repeller dual to $A_j$. 
\end{thm}
\proof We set $A_0:=\emptyset$ and $A_k:=\{x\in S\ |\ \mbox{there exists a solution } \sigma:\ZZ\to S\mbox{ through }x\mbox{ with }\alpha(\sigma)\subset M_1\cup M_2\cup\dots \cup M_{k}\}$, for $k\in\{1,2,\dots,n\}$.

The inclusions $\emptyset =A_0\subset A_1\subset A_2\subset\dots\subset A_n=S$ are straightforward.

We shall prove that the sets $A_k$ are attractors. We proceed inductively. First note that directly by the definition we have $A_n=S$; hence, our claim holds true for $k=n$. Now suppose that, for given $k\in\{0,1,\dots,n-1\}$, the set  $A_{k+1}$ is an attractor. We need to verify that $A_k$ is an attractor. 

We begin by showing that $\cl A_k$ and $M_{k+1}$ are disjoint. For contradiction, suppose that there exists $y\in \cl A_k\cap M_{k+1}$, and consider a sequence $\{y_t\}$ in $A_k$ convergent to $y$. For each $t\in\NN$ let $\sigma_t:\ZZ\to S$ be a solution for $F$ through $y_t$ with $\alpha(\sigma_t)\subset M_1\cup\dots\cup M_{k}$. Since $A_{k+1}$ is an attractor, $\sigma_t$ is a solution in $A_{k+1}$. We construct a solution $\sigma:\ZZ\to A_{k+1}$ through $y$. 
Fix $m\in\NN$, choose an increasing sequence $\{t_p\}\subset\NN$ such that $\sigma_{t_p}(l)$ is convergent in $A_{k+1}$ for each $l\in[-m,m]$, and set $\sigma^m(l):=\lim_{p\to\infty}\sigma_{t_p}(l)$. Clearly, $\sigma^m(0)=y$. Moreover, $\sigma_{t_p}(l+1)\in F(\sigma_{t_p}(l))$  for any $t_p\in\NN$ and $l,l+1\in[-m,m]$. Therefore, by the closed graph property of $F$ we obtain $\sigma^m(l+1)\in F(\sigma(l))$. This shows that $\sigma^m:[-m,m]\to A_{k+1}$ is a solution for $F$ through $y$. Proceeding inductively and using pointwise limits of solutions $\sigma_t$, we can extend $\sigma^m$ to a full solution $\sigma:\ZZ\to A_{k+1}$ for $F$ through $y$ (cf. \cite[Lemma 3.4]{BMMP2020}).   

Now, fix an arbitrary $l\in\{1,2,\dots,k\}$. Note that $F(M_l)\cap M_{k+1}=\emptyset$, because otherwise by the invariance of $M_l$ and $M_{k+1}$ we would have a solution $\tau:\ZZ\to S$ with $\alpha(\tau)\subset M_l$ and $\omega(\tau)\subset M_{k+1}$, which is in contradiction with the ordering of the Morse sets in $\cM$. Clearly, $M_l$ and $M_{k+1}$ are compact. The set $F(M_l)$ is compact as an image of a compact set under an  upper semicontinuous map. Therefore, by the upper semicontinuity of $F$ we can take compact and disjoint neighborhoods $N_l$ and $N' _{k+1}$ of $M_l$ and $M_{k+1}$, respectively, such that $F(N_l)\cap N' _{k+1}=\emptyset$. It follows that there exists a compact neighborhood $N_{k+1}$ of $M_{k+1}$ such that
$$
N_l\cap N_{k+1}=\emptyset\mbox{ and }F(N_l)\cap N_{k+1}=\emptyset
$$ 
for all $l\in\{1,2,\dots,k\}$. We set
$$
V_{k+1}:=A_{k+1}\setminus\bigcup _{i=1} ^{k+1} \Int N_i
$$
and
$$
V_{k}:=A_{k+1}\setminus\bigcup _{i=1} ^k \Int N_i.
$$
Since $\sigma_n(0)$ converges to $y\in M_{k+1}\subset\Int N_{k+1}$, we may assume that $\sigma_n(0)\in N_{k+1}$ for $n\in\NN$. Fix an $n\in\NN$. Observe that there exists a smallest $k_n\in\NN$ such that $\sigma_n(-k_n)\in V_{k+1}$, because $\alpha(\sigma_n)\subset M_1\cup\cdots\cup M_k$ and $F(N_1\cup\cdots\cup N_k)\cap N_{k+1}=\emptyset$.

{\em Case 1.} Sequence $\{k_n\}$ is bounded. Then we can find an $l\in\NN$ such that $\sigma_n(-l)\in V_{k+1}$ for all but finitely many $n\in\NN$. Passing in $\sigma_n(-l)$ with $n$ to infinity we obtain $\sigma(-l)\in V_{k+1}$, as $V_{k+1}$ is compact.  Since $\sigma(0)\in M_{k+1}$ and $M_{k+1}$ is invariant, there exists a solution $\tau:\ZZ\to M_{k+1}$ through $\sigma(0)$. Now, the concatenation
$$
\sigma '(m):=\left\{\begin{array}{rl}
\tau(m),&m\geq 0\\
\sigma(m),&m\leq 0
\end{array}\right.
$$
is a solution for $F$ in $A_{k+1}$ with $\omega(\sigma ')\subset M_{k+1}$ and $\alpha(\sigma')\subset M_1\cup\cdots\cup M_{k+1}$. However, $\sigma'(-l)=\sigma(-l)\notin M_{k+1}$, thus we cannot have $\alpha(\sigma')\subset M_{k+1}$. It follows that $\alpha(\sigma')\subset M_j$ for some $j\in\{1,2,\dots,k\}$, a contradiction.

{\em Case 2.} Sequence $\{k_n\}$ is unbounded. Without loss of generality we may assume that $k_n\geq n$. Observe that 
$$
\dom(F_{{V_k},{k_n}})\cap V_{k+1}\cap A_{k+1}\neq\emptyset
$$
for all $n\in\NN$. Taking into account that $\{\dom(F_{{V_k},{n}})\}$ is a decreasing sequence of compact sets, and $V_{k+1}$ and $A_{k+1}$ are compact, we get 
$$
\bigcap\left\{\dom(F_{{V_k},{k_n}})\cap V_{k+1}\cap A_{k+1}\,|\, n\in\NN\right\}\neq\emptyset.
$$
By the identity $\Inv^+(V_k,F)=\bigcap\left\{\dom(F_{{V_k},{n}})\,|\,n\in\ZZ ^+\right\}$ (cf. e.g. \cite[Lemma 4.9]{BM2016}, \cite[Lemma 2.8]{KM95}) we obtain 
$$
\Inv^+(V_k,F)\cap V_{k+1}\cap A_{k+1}\neq\emptyset.
$$ 
Take an $x\in \Inv^+(V_k,F)\cap V_{k+1}\cap A_{k+1}$ and a solution $\tau':\ZZ^+\to V_k$ for $F$ through $x$. Then $\omega(\tau')\subset V_{k}$, because $V_{k}$ is compact. This, in turn, means that $\omega(\tau')\subset M_{k+1}$, as $V_k\cap(M_1\cup M_2\cup\cdots\cup M_k)=\emptyset$ and $V_k\subset A_{k+1}$. Now, by the inclusion $x\in A_{k+1}$ and the invariance of $A_{k+1}$, we can extend $\tau'$ to the full solution $\tau:\ZZ\to A_{k+1}$. We have $\omega(\tau)\subset M_{k+1}$ and $\alpha (\tau)\subset M_1\cup M_2\cup\cdots\cup M_{k+1}$. Moreover, $\im\tau\not\subset M_{k+1}$, as $x\in V_{k+1}$. Therefore, we cannot have $\alpha(\tau)\subset M_{k+1}$, showing that $\alpha(\tau)\subset M_j$ for some $j\in\{1,2\dots,k\}$, a contradiction. 

This completes the proof that $\cl A_k$ and $M_{k+1}$ are disjoint.

Take a compact neighborhood $W_k$ of $A_k$ such that $W_k\cap M_{k+1}=\emptyset$. We claim that $W_k$ is an isolating neighborhood of $A_k$. Indeed, it suffices to verify the identity $A_k=\Inv(W_k, F)$. To this end consider $x\in\Inv(W_k,F)$ and $\sigma:\ZZ\to W_k$, a solution through $x$. By compactness of $W_k$ we have $\alpha(\sigma)\subset W_k$, and consequently $\alpha(\sigma)\subset M_1\cup\cdots\cup M_k$, showing that $x\in A_k$. We need to verify the other inclusion. It is straightforward to observe that $A_k$ is invariant. Thus, $A_k=\Inv(A_k,F)\subset \Inv(W_k, F)$. This completes the proof that $W_k$ is an isolating neighborhood for $A_k$. In particular, it follows that $A_k$ is compact. Moreover, using the same reasoning, one can easily see that $\Inv^-(W_k,F)=A_k$. 

Observe that $F(A_k)=A_k\subset\Int W_k$ and consider an open neighborhood $V$ of $A_k$ such that 
\begin{equation}\label{eq:v}
F(V)\subset\Int W_k.
\end{equation}
Using \cite[Lemma 2.11]{KM95} we can pick a compact neighborhood $A$ of $A_k$ such that $T_k:=F^+ _{W_k}(A) \subset V$. We shall show that $T_k$ is a trapping region for $A_k$. First note that the set $T_k$ is compact, because $\Inv^-(W_k,F)=A_k\subset A\subset W_k$ (cf. \cite[Lemma 2.9]{KM95}). Moreover, $A_k\subset\Int A\subset\Int T_k$. The identity $\Inv(T_k,F)=A_k$ is easily seen. Thus, $T_k$ is an isolating neighborhood for $A_k$. There remains to verify that $F(T_k)\subset T_k$.  Note that, directly by the definition, $T_k$ is positively invariant in $W_k$. Moreover, by (\ref{eq:v}) we have $F(T_k)=F(F^+ _{W_k}(A))\subset F(V)\subset W_k$. Consequently, $F(T_k)=F(T_k)\cap W_k\subset T_k$.

This completes the proof that $A_k$ is an attractor.

We will show the identity $M_j=A_j\cap A^* _{j-1}$ for $j\in\{1,\dots,n\}$. Fix $j\in\{1,\dots,n\}$. Let $x\in M_j$ and let $\sigma:\ZZ\to M_j$ be a solution through $x$. The set $M_j$ is compact, therefore $\alpha(\sigma)\subset M_j$. This means that $x\in A_j$. We shall verify that $x\in  A^* _{j-1}$. For contradiction suppose that $x\notin A^*_{j-1}$. Then, for any solution $\tau:\ZZ\to S$ we have $\omega(\tau)\cap A_{j-1}\neq\emptyset$. In particular, $\omega(\sigma)\cap A_{j-1}\neq\emptyset$, which in turn implies  $\omega(\sigma)\subset M_k$ for some $k\in\{1,2,\dots,j-1\}$. However, $\sigma$ is a solution in $M_j$; hence, $\omega(\sigma)\subset M_k \cap M_j$, a contradiction. For the proof of the other inclusion take an $x\in A_j\cap A^* _{j-1}$. Then there exists a solution $\sigma':\ZZ\to S$ through $x$ with $\alpha(\sigma')\subset M_1\cup\cdots\cup M_j$. We also have a solution $\sigma'':\ZZ\to S$ such that $\omega (\sigma'')\cap (M_1\cup\cdots\cup M_{j-1})=\emptyset$. Consequently, there is $k\geq j$ such that $\omega (\sigma'')\subset M_k$. Since both $\sigma'$ and $\sigma''$ are solutions through $x$, we can take their concatenation $\sigma:\ZZ\to S$ so that $\alpha(\sigma)\subset M_1\cup\cdots\cup M_j$ and $\omega (\sigma)\subset M_k$. Thus, by the ordering of the Morse sets $\alpha(\sigma)\subset M_j$ and $\omega (\sigma)\subset M_j$. It follows that $\im\sigma\subset M_j$, in particular $x\in M_j$. 

This completes the proof. 
\qed
\begin{prop}\label{prop:repattr-in-a}
Let $\cM:=\{M_1,M_2,\dots,M_n\}$  be a Morse decomposition of an isolated invariant set $S$ with respect to an upper semicontinuous $F:X\mto X$, and let  
$\emptyset =A_0\subset A_1\subset A_2\subset\dots\subset A_n=S$ be the sequence  of associated attractors. Then $(M_j,A_{j-1})$ is a repeller--attractor pair in $A_j$ for $j\in\{1,2,\dots,n\}$. 
\end{prop} 
\proof
Let $j\in\{1,2,\dots,n\}$ and let $T_{j-1}$ be a trapping region for $A_{j-1}$ in $S$. One can observe that  $T_{j-1}\cap A_j$ is a trapping region for $A_{j-1}$ in $A_j$, showing that $A_{j-1}$ is an attractor in $A_j$. Now, the assertion of the proposition is straightforward according to the identity $M_j=A_j\cap A^* _{j-1}$.
\qed

The following theorem may be viewed as a converse of Theorem \ref{thm:Morse_attr-rep}.
\begin{thm}\label{thm:Attr-rep-Morse}
Let $S$ be an isolated invariant set with respect to an upper semicontinuous map $F:X\mto X$, and let $\emptyset =A_0\subset A_1\subset A_2\subset\dots\subset A_n=S$ be a sequence of attractors in $S$. Then
$\cM:=\{M_1,M_2,\dots,M_n\}$, where  $M_j=A_j\cap A^* _{j-1}$ and $A^* _j$ is a repeller dual to $A_j$, for $j\in\{1,2,\dots,n\}$, is a Morse decomposition of $S$. Moreover, $A_k=\{x\in S\ |\ \mbox{there exists a solution } \sigma:\ZZ\to S\mbox{ through }x\mbox{ with }\alpha(\sigma)\subset M_1\cup M_2\cup\dots \cup M_{k}\}$, $k\in\{1,2,\dots,n\}$.
\end{thm}
\proof The idea of the proof is similar to its single valued counterpart (cf. e.g. \cite{RZ85}), but the details differ. We present the entire proof for the sake of completeness.  

We shall verify all the assertions of Definition \ref{def:Morse_dec}. Fix $i<j$. We have
$$
\begin{array}{rcl}
M_i\cap M_j&=&A_i\cap A^* _{i-1}\cap A_j\cap A^* _{j-1}\\
&\subset&A_i\cap A^* _{j-1}\subset A_{j-1}\cap A^* _{j-1}=\emptyset.
\end{array}
$$
Let $j\in\{1,2,\dots,n\}$ be fixed. We shall verify that $M_j$ is isolated and invariant with respect to $F$. Let $T_j$ be a trapping region for $A_j$. Recall that $N_{A^* _{j-1}}:=S\setminus\Int T_j$ is an isolating neighborhood for $A_{j-1}^*$ (cf. Lemma \ref{lem:rep-attr}). We claim that $N_j:=T_j\cap N_{A^* _{j-1}}$ is an isolating neighborhood for $M_j$. We have $M_j=A_j\cap A^* _{j-1}\subset A_j\subset \Int T_j$ and $M_j=A_j\cap A^* _{j-1}\subset A^* _{j-1}\subset \Int N_{A^* _{j-1}}$. Thus, $M_j\subset \Int T_j\cap \Int N_{A^* _{j-1}}=\Int  (T_j\cap N_{A^* _{j-1}})=\Int N_j$. There remains to show that $M_j=\Inv (N_j,F)$. To this end, consider an $x\in M_j$. Since $x\in A_j$, we can take a solution $\eta:\ZZ\to A_j$ with $\eta(0)=x$. Clearly, $\alpha (\eta)\subset A_j$. On the other hand $x\in A^* _{j-1}$, therefore there exists a solution $\tau:\ZZ\to S$ with $\tau(0)=x$ such that $\omega(\tau)\cap A_{j-1}=\emptyset$. Define $\sigma:\ZZ\to S$ by
$$
\sigma(k):=\left\{\begin{array}{rl}
\eta(k),&k\leq 0,\\
\tau(k), &k\geq 0.
\end{array}\right.
$$
One easily sees that $\sigma$ is a solution with respect to $F$ through $x$. Note that $\omega(\sigma)=\omega(\tau)$; hence, $\omega(\sigma)\cap A_{j-1}=\emptyset$. Then, by Lemma \ref{lem:omega_in_A}, $\im\sigma\cap T_{A _{j-1}}=\emptyset$  and by Lemma \ref{lem:rep-attr} the inclusion $\im\sigma\subset A^* _{j-1}$ follows. Moreover $\alpha(\sigma)=\alpha(\eta)\subset A_j$, and by Lemma \ref{lem:rep-attr-sol}(iv), $\im\sigma\subset A_j$. Consequently, $\im\sigma\subset A_j\cap A^* _{j-1}=M_j$ showing that $M_j$ is invariant with respect to $F$. Therefore, $M_j=\Inv(M_j, F)\subset \Inv (N_j,F)$. 
For the proof of the other inclusion consider an $x\in \Inv(N_j,F)$ and a solution $\sigma:\ZZ\to N_j$ through $x$. We have $\im\sigma\subset N_j\subset T_j$, and in turn $\im\sigma\subset A_j$, as $T_j$ is a trapping region for $A_j$. Similarly, $\im\sigma\subset N_j\subset N_{A^* _{j-1}}$ and we have the inclusion $\im\sigma\subset A^* _{j-1}$, because $N_{A^* _{j-1}}$ is an isolating neighborhood for $A^* _{j-1}$. As a consequence, $\im\sigma\subset A_j\cap A^* _{j-1}=M_j$. In particular $x\in M_j$. We have proved that $M_j$ is an isolated invariant set, which also justifies its compactness. The proof of Definition \ref{def:Morse_dec}(1) is complete.

In order to prove Definition \ref{def:Morse_dec}(2) consider a solution $\sigma:\ZZ\to S$. Since the sequence of attractors is increasing and $A_n=S$, there exists a smallest positive integer $i$ such that $\omega(\sigma)\subset A_i$. Similarly, there is a largest $j\in\NN$, $j<n$, with $\alpha(\sigma)\subset A^* _j$, as the sequence of the dual repellers is decreasing and $A^* _0=S$. Then, $\omega(\sigma)\not\subset A_{i-1}$, which according to Lemma \ref{lem:omega_in_A} implies $\im\sigma\cap T_{i-1}=\emptyset$. Hence, by Lemma \ref{lem:rep-attr} we get 
\begin{equation}\label{eq:im_a*}
\im\sigma\subset A^* _{i-1}.
\end{equation}
In particular, $\omega(\sigma)\subset A^* _{i-1}$. This, along with the inclusion $\omega(\sigma)\subset A_i$ implies 
\begin{equation}\label{eq:omega_mi}
\omega(\sigma)\subset M_i.
\end{equation}

By the choice of $j$ we have $\alpha(\sigma)\subset A^* _j$ and $\alpha(\sigma)\not\subset A^* _{j+1}$. We claim that 
\begin{equation}\label{eq:im_aj+1}
\im\sigma\subset A_{j+1}.
\end{equation} 
Indeed, otherwise we would have $\im\sigma\not\subset T_{A_{j+1}}$ and, by Lemma \ref{lem:rep-attr-sol}(i), $\alpha(\sigma)\subset A^* _{j+1}$, a contradiction.

Observe that $j+1\geq i$. If not, then by (\ref{eq:im_aj+1}) $\im\sigma\subset A_{j+1}\subset A_{i-1}$, which along with (\ref{eq:im_a*}) yields $\im\sigma\subset A_{i-1}\cap A^* _{i-1}=\emptyset$, a contradiction.

If $j+1=i$ then by (\ref{eq:im_aj+1}) we have $\im\sigma\subset A_{j+1}= A_{i}$. Now, using (\ref{eq:im_a*}) we obtain $\im\sigma\subset A_{i}\cap A^* _{i-1}=M_i$.

If $j+1>i$ then we have $\alpha(\sigma)\subset A^* _{j}$ and, by (\ref{eq:im_aj+1}), $\alpha(\sigma)\subset A_{j+1}$. Thus, $\alpha(\sigma)\subset A_{j+1}\cap A^* _{j}=M_{j+1}$. This along with (\ref{eq:omega_mi}) completes the proof of Definition \ref{def:Morse_dec}(2).

For the proof of Definition \ref{def:Morse_dec}(3) fix a $j\in\{1,2,\dots,n\}$ and take a solution $\sigma:\ZZ\to S$ with $\alpha(\sigma)\cup\omega(\sigma)\subset M_j$. Then the inclusion $\im\sigma\subset A_j$ follows from Lemma \ref{lem:rep-attr-sol}(iv), as $\alpha(\sigma)\subset A_j$. Similarly, by Lemma \ref{lem:rep-attr-sol}(iii) and  the inclusion $\omega(\sigma)\subset A^* _{j-1}$ we get $\im\sigma\subset A^* _{j-1}$. Consequently, $\im\sigma\subset M_j$. 

We have proved that the family $M$ is a Morse decomposition of $S$.

There remains to justify the last assertion of the theorem. Fix a $k\in\{1,2,\dots,n\}$. Take an $x\in A_k$ and a solution $\sigma:\ZZ\to A_k$ through $x$. Clearly, $\alpha(\sigma)\subset A_k$. Let $i\leq k$ be the smallest integer such that $\alpha(\sigma)\subset A_i$. Then $\alpha(\sigma)\not\subset A_{i-1}$. This means that $\im\sigma\not\subset T_{i-1}$, which along with Lemma \ref{lem:rep-attr-sol}(i) yields $\alpha(\sigma)\subset A^*_{i-1}$. Thus, $\alpha(\sigma)\subset A_i\cap A^* _{i-1}=M_i$. For the proof of the other inclusion consider an $x\in S$ and a solution $\sigma:\ZZ\to S$ through $x$ with $\alpha(\sigma)\subset  M_1\cup M_2\cup\dots \cup M_{k}$. There is $i\leq k$ such that $\alpha(\sigma)\subset  M_i$. Then $\alpha(\sigma)\subset  A_i$, which along with Lemma \ref{lem:rep-attr-sol}(iv) shows that $\im\sigma\subset A_i$. In particular, $x\in A_i\subset A_k$.

This completes the proof.
\qed
\section{Index triples}\label{sec:triples}   
Throughout this section assume that $X$ is a locally compact metric space and $F:X\mto X$ is an upper semicontinuous multivalued map determined by a morphism.

\begin{lem}\label{lm:l1}
Let $N$ be an isolating neighborhood for $F$ and let a pair $P=(P_1,P_2)$ of compact sets $P_2\subset P_1\subset N$ satisfy conditions (a) and (d) of Definition \ref{defn:wip}. Then $P$ satisfies condition (b).
\end{lem}
\proof First we show that $\bd_F P_1\subset \bd N$. For contradiction assume that there exists a $y\in \bd_F P_1\setminus \bd N$. Then $y\in\cl(F(P_1)\setminus P_1)$, and we can consider a sequence $\{y_n\}\subset F(P_1)\setminus P_1$ convergent to $y$. Observe that $y\in\Int N$, because $y\in P_1\subset N$ and $y\notin \bd N$. Therefore, $y_n\in\Int N$ for large enough $n\in\NN$. This along with $y_n\in F(P_1)$ and (a) implies $y_n\in P_1$, a contradiction. 

To prove inclusion (b) assume the contrary and consider an $x\in  \bd_F P_1\setminus P_2$. Then $x\in P_1\setminus P_2\subset \Int N$, by property (d). On the other hand, $x\in \bd_F P_1\subset \bd N$, a contradiction.  
\qed

The next lemma is straightforward.
\begin{lem}\label{lm:pos_inv}
Let $A,B,N,M$ be subsets of $X$. The following holds true:
\begin{itemize}
\item[(i)] if $A,B$ are positively invariant with respect to $F$ in $N$ then so are $A\cup B$ and $A\cap B$,
\item[(ii)] if $A$ and $B$ are positively invariant with respect to $F$ in $N$ and $M$, respectively, then $A\cap B$ is positively invariant with respect to $F$ in $N\cap M$,
\item[(iii)] if $A$  is positively invariant with respect to $F$ in $N$ and $B$ is positively invariant with respect to $F$, then $A\cap B$ is positively invariant with respect to $F$ in $N$,
\item[(iv)] if $A$ is positively invariant with respect to $F$ in $M$ and $N\subset M$ then $A$ is positively invariant with respect to $F$ in $N$.
\end{itemize}
\end{lem}
\begin{defn}\label{defn:F-pair} 
{\rm We will say that a pair $R=(R_1,R_2)$ of compact sets is an {\em $F$--pair} if there is a compact set $M$ with $R_2\subset R_1\subset M$ and 
\begin{itemize}
\item[(Fp1)] $R_1$, $R_2$ are positively invariant with respect to $F$ in $M$,
\item[(Fp2)] $\cl (R_1\setminus R_2)$ is an isolating neighborhood,
\item[(Fp3)] $R_1\setminus R_2\subset \Int M$.
\end{itemize}
}
\end{defn}
Note that a weak index pair is a special $F$--pair.

Given an $F$--pair $R$ in $M$ we set
$$
T_M(R):=(T_{M,1}(R),T_{M,2}(R)):=(R_1\cup (X\setminus\Int M),R_2\cup (X\setminus\Int M)).
$$
\begin{lem}\label{lem:F-pair}
If $R$ is an $F$--pair in $M$, then 
\begin{itemize}
\item[(i)] $F(R)\subset T_M(R)$,
\item[(ii)] the inclusion $i_{R}:=i_{R,T_M(R)}:R\to T_M(R)$ induces an isomorphism in the Alexander-Spanier cohomology.
\end{itemize}
\end{lem}
\proof Property (i) follows from the positive invariance of $R_1$ and $R_2$ with respect to $F$ in $M$. Since $T_{M,1}(R)\setminus T_{M,2}(R)=R_1\setminus R_2$, inclusion $i_{R}$ is an excision, and property (ii) follows.
\qed

By Lemma \ref{lem:F-pair} we can define an endomorphism $I_R:H^*(R)\to H^*(R)$ by
$$
I_R:=F_R ^*\circ (i_R ^*)^{-1},
$$ 
where $F_R$ stands for the restriction of $F$ to the domain $R$ and the codomain $T_M(R)$. In the following endomorphism $I_R$ will be called an {\em index map} associated with an $F$--pair $R$.

Note that if an $F$--pair $R$ is a weak index pair then the above notion of the index map coincides with that for a weak index pair.
\begin{prop}\label{prop:F-pair_index}
Let $R$ be an $F$--pair in $M$, and let $S:=\Inv(\cl(R_1\setminus R_2),F)$. If $N\subset M$ is an isolating neighborhood of $S$ such that $R_1\setminus R_2\subset\Int N$ then 
\begin{itemize}
\item[(i)] $P:=R\cap N$ is a weak index pair for $F$ in $N$,
\item[(ii)] the inclusion $i_{PR}:P\to R$ induces an isomorphism in the Alexander-Spanier cohomology,
\item[(iii)] index maps $I_{F_P}$ and $I_R$ are conjugate.
\end{itemize}
As a consequence, $C(S,F)=L(H^*(R),I_R)$.
\end{prop}
\proof 
We shall verify that $P$ is a weak index pair. First observe that, by the positive invariance of $R$ in $M$, the inclusion $N\subset M$, and Lemma \ref{lm:pos_inv}(iv), $P$ is positively invariant in $N$, i.e. $P$ satisfies property (a) of Definition \ref{defn:wip}. Properties (c) and (d) are straightforward. Now, property (b) is justified by Lemma \ref{lm:l1}.

Since $R_1\setminus R_2\subset\Int N$, we have $P_1\setminus P_2=(R_1\setminus R_2)\cap N=R_1\setminus R_2$. Therefore, inclusion $i_{PR}$ induces an isomorphism, as an excision. 

There remains to verify (iii). Consider the commutative diagram
$$
\begin{tikzcd}
(P_1,P_2)\ar{d}{i_{PR}}\ar{r}{F_P}&(T_{N,1}(P),T_{N,2}(P))&(P_1,P_2)\ar{d}{i_{PR}}\ar{l}[swap]{i_P}
\\[3ex]
(R_1,R_2)\ar{r}{F_{R}} & (T_{M,1}(R),T_{M,2}(R))\ar{u}{j} & (R_1,R_2)\ar{l}[swap]{i_{R}}
\end{tikzcd}
$$
in which $i_P$, $i_{R}$, $i_{PR}$ and $j$ are inclusions. Recall that $i_P$ and $i_{R}$ and $i_{PR}$ induce isomorphisms in cohomology by the strong excision property. And, so does $j$. The commutativity of the diagram shows that the index maps $I_{F_P}$ and $I_R$ are conjugate. This completes the proof.
\qed

\begin{lem}\label{lm:l2}
Assume $S_1\subset S_2$ are isolated invariant sets with respect to $F$, with isolating neighborhoods $N_1$ and $N_2$, respectively. Then, $N_1\cap N_2$ is an isolating neighborhood for $S_1$.
\end{lem}
\proof
We have $S_1=\Inv (S_1,F)\subset \Inv (N_1\cap N_2,F)$, because $S_1$ is invariant and $S_1\subset N_1\cap N_2$. On the other hand, $\Inv(N_1\cap N_2,F)\subset\Inv(N_1,F)=S_1$. Thus, $S_1=\Inv(N_1\cap N_2,F)$. There remains to verify that $S_1\subset\Int(N_1\cap N_2)$. To this end observe that $S_1\subset\Int N_1$ and $S_1\subset S_2\subset \Int N_2$. Therefore, $S_1=S_1\cap S_2\subset \Int N_1\cap\Int N_2=\Int(N_1\cap N_2)$.
\qed
\begin{thm}\label{thm:index-triple}
Let $S$ be an isolated invariant set with respect to $F$ and let $N$ be its isolating neighborhood. Assume that $(A^*,A)$ is a repeller--attractor pair in $S$. Then, there exist a triple $(P_0,P_1,P_2)$ of compact subsets $P_2\subset P_1\subset P_0$ of $N$ such that
\begin{itemize}
\item[(i)] $(P_0,P_2)$ is a weak index pair for $F$ and $C(S,F)=L(H^*(P_0,P_2),I_{(P_0,P_2)})$,
\item[(ii)] $(P_1,P_2)$ is an $F$--pair and $C(A,F)=L(H^*(P_1,P_2),I_{(P_1,P_2)})$,
\item[(iii)] $(P_0,P_1)$ is an $F$--pair and $C(A^*,F)=L(H^*(P_0,P_1),I_{(P_0,P_1)})$.
\end{itemize}
\end{thm}
\proof
Let $(P_0,P_2 ')$ be a weak index pair for $S$ and $F$ in an isolating neighborhood $N$. Consider a trapping region $M'$ for attractor $A$. By Lemma \ref{lm:l2} the set $M:=M'\cap N$ is an isolating neighborhood of $A$. So, we can take a weak index pair $Q'=(Q_1',Q_2 ')$ for $A$ and $F$ in $M$. Set $Q_i:=Q_i'\cap P_0$, $i\in\{1,2\}$, and $Q:=(Q_1,Q_2)$.

First we prove that $Q$ is a weak index pair for $A$ in $M$. Fix $i\in\{1,2\}$. Since $Q_i'$ is positively invariant in $M$ and $P_0$ is positively invariant in $N$, by Lemma \ref{lm:pos_inv}(ii) $Q_i=Q_i'\cap P_0$ is positively invariant in $M\cap N=M$. This shows that $Q$ satisfies property (a) of Definition \ref{defn:wip}. Observe that $A\subset S\subset \Int P_0$ and $A\subset\Int(Q_1'\setminus Q_2')$, as $(P_0,P_2')$ and $(Q_1,Q_2')$ are weak index pairs for $S$ and $A$, respectively. Therefore, $A\subset \Int(Q_1'\setminus Q_2')\cap\Int P_0=\Int((Q_1'\setminus Q_2')\cap P_0)=\Int(Q_1\setminus Q_2)$, showing (c). Property (d) follows from the straightforward inclusions $Q_1\setminus Q_2=(Q_1'\cap P_0)\setminus (Q_2'\cap P_0)=(Q_1'\setminus Q_2')\cap P_0\subset Q_1'\setminus Q_2'\subset\Int M$. Now, property (b) is a consequence of (a), (d), and Lemma \ref{lm:l1}.

Set
\begin{equation}\label{eq:def_p1p2}
P_1:=Q_1\cup P_2'\mbox{ and }P_2:=Q_2\cup P_2'.
\end{equation}
 We shall show that $(P_1, P_2)$ is an $F$--pair in $N$. Clearly, $P_2\subset P_1$ are compact subsets of $N$. Property (Fp3) in Definition \ref{defn:F-pair} is straightforward, because $P_1\setminus P_2\subset Q_1\setminus Q_2$ and $Q$ is a weak index pair for $F$ in $M\subset N$. For the proof of (Fp2) observe that $\cl(P_1\setminus P_2)\subset \cl(Q_1\setminus Q_2)\subset M$, therefore $\Inv(\cl(P_1\setminus P_2),F)\subset \Inv(M,F)=A$. Moreover, $A\subset\Int(Q_1\setminus Q_2)$ and $A\subset\Int(P_0\setminus P'_2)$, as $Q$ is a weak index pair for $A$, $(P_0,P'_2)$ is a weak index pair for $S$, and $A\subset S$. Thus, $\Inv(\cl(P_1\setminus P_2),F)\subset\Int(Q_1\setminus Q_2)\cap\Int(P_0\setminus P'_2)=\Int(P_1\setminus P_2)$. There remains to verify (Fp1). For $i\in\{1,2\}$ we have $F(P_i)\cap N=F(Q_i\cup P_2')\cap N=((F(Q_i)\cap N)\cup (F(P_2'))\cap N)$. By the positive invariance of $P_2'$ in $N$ we have the inclusion $F(P_2')\cap N\subset P_2'$. Recall that $M'$ is an attracting neighborhood and $Q_i\subset M'$; hence, $F(Q_i)\subset M'$. Taking this into account we obtain $F(Q_i)\cap N\subset F(Q_i)\cap N\cap M'=F(Q_i)\cap M\subset Q_i$, where the last inclusion is a consequence of the positive invariance of $Q_i$ in $M$. Finally, we have  $F(P_i)\cap N\subset Q_i\cup P_2'=P_i$, i.e. $P$ satisfies (Fp1). This shows that $(P_1,P_2)$ is an $F$--pair in $N$.

Notice that $\Inv(P_1\setminus P_2,F)=A$ and  $M$ is an isolating neighborhood of $A$ with $P_1\setminus P_2\subset \Int M$. Therefore, by Proposition \ref{prop:F-pair_index} we infer that $C(A,F)=L(H^*(P_1,P_2),I_{(P_1,P_2)})$. This completes the proof of (ii).

We shall show that $P:=(P_0, P_2)$ is a weak index pair for $S$ and $F$ in $N$. Recall that $(P_0,P_2')$ is a weak index pair for $F$ in $N$. Thus, in order to verify that $P$ satisfies property (a) in Definition \ref{defn:wip} it suffices to show that $P_2=Q_2\cup P_2'$ is positively invariant in $N$. Since the sets $Q_2$ and $P_2'$ are positively invariant in $N$, our claim follows from Lemma \ref{lm:pos_inv}(i). We have $P_1\setminus P_2=P_1\setminus (Q_2\cup P_2')\subset P_1\setminus P_2'\subset\Int N$, showing property (d). Property (b) is a consequence of (a), (d) and Lemma \ref{lm:l1}. There remains to verify property (c), i.e. the inclusion $S\subset\Int (P_1\setminus P_2)$. Since $S\subset\Int (P_1\setminus P_2')$, for contradiction suppose that there is an $x\in S\cap Q_2$. Consider a solution $\sigma$ through $x$ in $S$. Since $Q_2$ is compact and positively invariant in $N$, we infer that $\emptyset\neq\omega (\sigma)\subset Q_2$. Moreover, $Q_2\subset M$ and $\omega (\sigma)$ is invariant with respect to $F$, hence we have $\omega(\sigma)=\Inv(\omega(\sigma),F)\subset \Inv(M,F)=A$, showing that $A\cap Q_2\neq\emptyset$. However, $(Q_1,Q_2)$ is a weak index pair for $A$, a contradiction. 

Since $(P_0,P_2)$ is a weak index pair for $S$ and $F$ in $N$, the identity $C(S,F)=L(H^*(P_0,P_2),I_{(P_0,P_2)})$ is straightforward. This completes the proof of (i).

Now we focus on a repeller $A^*$ dual to $A$ in $S$. We shall verify that $(P_0,P_1)$ is an $F$--pair in $N$. Clearly, the sets $P_1\subset P_0$ are compact subsets of $N$. Recall that we have already justified the positive invariance of $P_0$ and $P_1$ in $N$, that is property (Fp1). Moreover, by the inclusions $P_0\setminus P_1\subset P_0\setminus P_2\subset\Int N$, (Fp3) follows. 
We claim that $\cl(P_0\setminus P_1)$ is an isolating neighborhood for $A^*$. According to (ii) we have $A\subset\Int P_1$, showing that $\cl(P_0\setminus P_1)\cap A=\emptyset$. Thus, according to Lemma \ref{lem:rep_neighborhood}, it suffices to verify that 
\begin{equation}\label{eq:a*in}
A^*\subset \Int (P_0\setminus P_1).
\end{equation}
To this end consider an $x\in A^*$ and a solution $\sigma$ through $x$ such that $\omega(\sigma)\cap A=\emptyset$. Observe that $x\notin Q_1$, because otherwise we would have $\emptyset\neq\omega(\sigma)\subset Q_1$, and consequently $\omega(\sigma)\subset A$, leading a contradiction. Thus, by the inclusion $A^*\subset S\subset\Int (P_0\setminus P_2)$ and (\ref{eq:def_p1p2}), that is the definition of $P_1$, inclusion (\ref{eq:a*in}) follows. 

Now, using (\ref{eq:a*in}), (\ref{eq:def_p1p2}) and Lemma \ref{lem:rep_neighborhood}, one can easily prove that $K:=N\setminus\Int Q_1$ is an isolating neighborhood of $A^*$. We claim that $(P_0\cap K,P_1\cap K)$ is a weak index pair for $F$ in $K$. Indeed, in view of Proposition \ref{prop:F-pair_index}(i) it suffices to justify that $P_0\setminus P_1\subset \Int K$, which follows from the inclusions $P_0\setminus P_1\subset\Int N$ and $P_0\setminus P_1\subset P_0\setminus Q_1$.
 
Eventually, according to Proposition \ref{prop:F-pair_index} we have the identity\linebreak $C(A^*,F)=L(H^*(P_0,P_1),I_{(P_0,P_1)})$, which proves (iii).

This completes the proof.
\qed
\begin{ex}\label{ex:wip_triple}
{\rm Consider multivalued map given in Example \ref{ex:F} (see Figure \ref{fig:Morse_dec}). We present some examples of weak index triples. 
\begin{itemize}
\item[(a)] Consider a repeller--attractor pair $(A_1 ^*, A_1)=([0,m]\cup[n,o]\cup [p,q]\cup[g,h],[d,e]\cup[j,k])$ in $S=[0,k]$ (Figure \ref{fig:Morse_dec}(C)). Note that $([0,l],\emptyset)$ is a weak index pair for $S$ in $[0,l]$ and  $([c,f]\cup[i,l],\emptyset)$ is a weak index pair for $A_1$ in $[c,f]\cup[i,l]$. Observe that $P_0:=[0,l]$, $P_1:=[c,f]\cup[i,l]$ and $P_2:=\emptyset$ is a triple satisfying the assertions of Theorem \ref{thm:index-triple}.
\item[(b)] Observe that for the repeller--attractor pair $(A_2 ^*,A_2)=([0,a],[d,k])$ in $S$ (Figure \ref{fig:Morse_dec}(D)), the weak index triple is given by $P_0:=[0,l]$, $P_1:=[c,l]$ and $P_2:=\emptyset$.
\item[(c)] Note that $(A_2 ^*, M_2)$ is a repeller--attractor pair in $A_1 ^*$ (see Figure \ref{fig:Morse_dec}(C) and (D). Then, $([0,c]\cup[f,i],\{c\}\cup\{f\}\cup\{i\})$ is a weak index pair for $A_1 ^*$, and $([f,i],\{f\}\cup\{i\})$ is a weak index pair for $M_2$. One can verify that the triple $([0,c]\cup[f,i],[f,i]\cup \{c\},\{c\}\cup\{f\}\cup\{i\})$ satisfies the assertions.
\end{itemize} 
}
\end{ex}
\section{Morse equation}\label{sec:Morse_eq}
The aim of this section is to prove the Morse equation. 

As in the preceding section, $X$ is a locally compact metric space and $F:X\mto X$ is an upper semicontinuous multivalued map determined by a morphism.
\begin{prop}\label{prop:Q}
Assume $(A^*,A)$ is a repeller--attractor pair in an isolated invariant set $S$ with respect to an upper semicontinuous map $F:X\mto X$. Then
$$
p(t,A^*)+p(t,A)=p(t,S)+(1+t)Q(t)
$$
for some power series $Q$ with nonnegative integer coefficients.
\end{prop}
\proof The proof runs along the lines of the proof of \cite[Theorem 5.5]{MMr91} with the use of Theorem \ref{thm:index-triple} instead of \cite[Theorem 4.3]{MMr91}.
\qed
\begin{thm}\label{thm:Morse_eq}
Let $\cM:=\{M_1,M_2,\dots,M_n\}$  be a Morse decomposition of an isolated invariant set $S$ with respect to an upper semicontinuous $F:X\mto X$. Then
\begin{equation}\label{eq:M_eq}
\sum _{i=1} ^n p(t,M_i)=p(t,S)+(1+t)Q(t),
\end{equation}
where $Q$ is a formal power series with nonnegative integer coefficients. Moreover, 
$$
Q(t)=\sum _{i=1} ^n Q_i(t),
$$
where 
\begin{equation}\label{eq:Q_i}
(1+t)Q_i(t)=p(t,M_i)+p(t,A_{i-1})-p(t,A_i).
\end{equation}
If $Q_i\neq 0$ then there exists a solution $\sigma:\ZZ\to S$ with $\alpha(\sigma)\subset M_i$ and $\omega(\sigma)\subset M_j$ for some $j<i$.
\end{thm} 
\proof
The proof runs similarly as the proof of \cite[Theorem 5.6]{MMr91}. By Theorem \ref{thm:Morse_attr-rep} and Proposition \ref{prop:repattr-in-a} we can consider the associated sequence of attractors such that $(M_i, A_{i-1})$ is a repeller--attractor pair in $A_i$. Applying Proposition \ref{prop:Q} to this pair we find a formal power series $Q_i$ with nonnegative integer coefficients satisfying  
\begin{equation}\label{eq:Qi}
p(t,M_i)+p(t,A_{i-1})=p(t,A_i)+(1+t)Q_i(t).
\end{equation}
Summing these equations over $i\in\{1,2,\dots,n\}$ and taking into account that $A_n=S$ and $A_0=\emptyset$ we obtain the first assertion.

In order to prove the second assertion, for contradictions suppose that there is no solution $\sigma:\ZZ\to S$ with $\alpha(\sigma)\subset M_i$ and $\omega(\sigma)\subset M_j$, where $j<i$. This means that $A_i$ is a disjoint union of $M_i$ and $A_{i-1}$,  $F(M_i)\cap A_{i-1}=\emptyset$, and $F(A_{i-1})\cap M_i=\emptyset$. By the additivity property of the Conley index \cite[Theorem 5.3]{B2017} we obtain 
$$
C(A_i,F)=C(M_i,F)\times C(A_{i-1},F).
$$
In particular,
$$
\dim C^q(A_i,F)=\dim C^q(M_i,F)+\dim C^q(A_{i-1},F).
$$
Now, there remains to multiply the above equalities by $t^q$ and sum over $q\in\NN$ in order to get
$$
p(t,A_i)=p(t,M_i)+p(t,A_{i-1}),
$$
which, along with (\ref{eq:Qi}), implies $Q_i= 0$, a contradiction. 
\qed
\begin{ex}\label{ex:M_eq}
{\rm Consider Morse decomposition of an isolated invariant set $S=[0,k]$ with respect to a multivalued map presented in Example \ref{ex:F} and Figure \ref{fig:Morse_dec}. One can verify that $p(t,M_1)=2$, $p(t,M_2)=t$, $p(t,M_3)=0$, and $p(t,S)=1$, so that the Morse equation has the form 
$$
2+t+0=1+(1+t)Q(t),
$$ 
showing that $Q(t)=1$. Moreover, $p(t,A_0)=0$, $p(t,A_1)=2$, $p(t,A_2)=1$ and $p(t,A_3)=1$. Therefore,  by (\ref{eq:Q_i}) we obtain that $Q_1(t)=Q_3(t)=0$ and $Q_2(t)=1$. Thus, the Morse equation shows that there exists a solution $\sigma:\ZZ\to S$ with $\alpha(\sigma)\subset M_2$ and $\omega(\sigma)\subset M_1$.
}
\end{ex}
Assume, we are given a Morse decomposition with respect to a flow on $X$. Recall that, according to \cite[Theorem 1]{MMr90}, any isolated invariant set with respect to a flow is isolated and invariant with respect to time-one-map. Note that this extends over the Morse decomposition. Thus, by Theorem \ref{thm:Morse_eq} we get Morse equation (\ref{eq:M_eq}) for the flow.   
Therefore, arguing as in \cite{RZ85} or \cite{RS99}, as a corollary we obtain the classical Morse inequalities.
\begin{cor}
Let $M$ be a $d$-dimensional compact manifold with $k$-th Betti number $\beta_k$, $k\in \{0,1,2\dots,d\}$. 
Let nondegenerate critical points $\{x_1,x_2,\dots,x_n\}$ with respect to a given gradient flow on $M$ constitute a Morse decomposition of $M$, and let $\lambda(x_j)$ stand for the Morse index of $x_j$.  
Then
$$
\sum _{j=1} ^n t^{\lambda(x_j)}=\sum _{j=0} ^d \beta_jt^j+(1+t)Q(t),
$$ 
where $Q$ is a formal power series with nonnegative integer coefficients. 
\end{cor}
\section*{Acknowledgements}
The paper was finalized during author's visit to the Rutgers University, USA. It is my pleasure to thank the staff of the Department of Mathematics, especially Professor Konstantin Mischaikow for the invitation and for inspiring scientific atmosphere.
  

\end{document}